\documentclass[a4paper,11pt]{amsart}
\usepackage{amssymb, amsmath,latexsym,amsfonts,amsbsy, amsthm,mathtools,graphicx,CJKutf8,CJKnumb,CJKulem,color,geometry}
\usepackage{verbatim}
\def\l{\left}
\def\r{\right}
\def\f{\frac}

\def\az{\alpha}
\def\lz{\lambda}
\def\Lz{\Lambda}
\def\az{\alpha}
\def\rz{\rho}
\def\ez{\epsilon}
\def\bz{\beta}
\def\dz{\delta}
\def\gz{\gamma}
\def\tz{\theta}
\def\sz{\sigma}

\def\beq{\begin{equation}}
\def\beq*{\begin{equation*}}

\author{Ovidiu Savin and Qian Zhang}
\title{\bf{Boundary H\"older gradient estimates for the Monge-Amp\`ere equation}}
\date{}

\theoremstyle{plain}\newtheorem{defn}{Definition}[section]
\theoremstyle{plain}\newtheorem{thm}{Theorem}[section]
\theoremstyle{plain}
\theoremstyle{plain}
\theoremstyle{plain}\newtheorem{lem}{Lemma}[section]
\theoremstyle{plain}

\numberwithin{equation}{section}

\begin{document}
\maketitle
\noindent {\bf{Abstract}.}\quad
We investigate global H\"older gradient estimates for solutions to the Monge-Amp\`ere equation
$$\mathrm{det}\;D^2 u=f\quad\mathrm{in}\;\Omega,$$
where the right-hand side $f$ is bounded away from $0$ and $\infty$. We consider two main situations when a) the domain $\Omega$ is uniformly convex and b) $\Omega$ is flat.

\bigskip

\section{Introduction}\label{s1}

In this paper, we consider boundary H\"older gradient estimate for solutions to the Dirichlet problem 
\begin{equation}\label{s1: eq 1}
\l\{\begin{array}{rclcl}
\mathrm{det}\;D^2 u&=&f&&\mathrm{in}\;\Omega,\\
u&=&\varphi&&\mathrm{on}\;\partial\Omega,
\end{array}
\r.
\end{equation}
where $\Omega$ is a convex domain in $\mathbb{R}^n$ and $0<\lz\le f\le\Lz$ for some constants $\lz,\Lz$.

The regularity of solutions for the Monge-Amp\`ere equation has been extensively studied by many authors, see for instance \cite{CY, P, W1, C1, C2, C3, H, C, CNS, I, K, S2, S4, S3, S5, TW} and references therein.

Concerning gradient H\"older estimates, Caffarelli proved in \cite{C1,C3} that solutions $u$ of \eqref{s1: eq 1} which are strictly convex satisfy $u \in C^{1,\delta}$ in the interior of $\Omega$, for some small $\delta>0$ depending on $\lambda$, $\Lambda$ and the dimension $n$. Moreover, the strict convexity of solutions can be guaranteed if the boundary data is above a critical  regularity level $\varphi \in C^{1,\beta}$ with $\beta>1-\f{2}{n}$. This exponent is optimal in view of Pogorelov's famous example of singular solutions in \cite{P}. However, as we will see later, even in this case the $C^{1,\delta}$ norm of $u$ may degenerate near the boundary of $\Omega$.

Here we investigate the $C^{1,\alpha}$ estimates up to the boundary of $\Omega$, under minimal conditions on the domain and the boundary data. While there is a rich literature addressing $C^{2,\alpha}$ boundary estimates for solutions of \eqref{s1: eq 1}, to the authors knowledge there is no work concerning sharp $C^{1,\alpha}$ boundary estimates which we discuss in this paper.  

We consider two main situations when a) the domain $\Omega$ is uniformly convex and b) $\Omega$ is flat.
In both cases we state two results similar in nature, one of them regarding the pointwise $C^{1,\alpha}$ estimate at a point on $\partial \Omega$ and the other one about the global version of this estimate.

For uniformly convex domains, Theorem \ref{s1: thm 1} below states that if $\partial\Omega$ and the boundary data $\varphi$ are pointwise $C^{2,\az}$ at a boundary point for some $\az\in(0,1)$, then the solution $u$ is $C^{1,\dz_0}$ at this point for some small $\dz_0>0$. 

\begin{thm}\label{s1: thm 1}
Let $u:\overline{\Omega}\to\mathbb{R}$ be a convex, continuous solution to \eqref{s1: eq 1}. Assume $\Omega\subset\mathbb{R}^n_+, 0\in\partial\Omega$, $\Omega$ is uniformly convex at $0$, and $\partial\Omega,\varphi\in C^{2,\az}(0)$; i.e., we assume that on $\partial\Omega$ we have
$$x_n=q(x')+O(|x'|^{2+\az}),$$
$$\varphi(x)=p(x')+O(|x'|^{2+\az}),$$
where $p(x'),q(x')$ are quadratic polynomials. Then 
$$u\in C^{1,\dz_0}(0)$$
for some constant $\dz_0>0$ depending only on $n$ and $\az$. 
\end{thm}

For the definition of $C^{1,\dz}(0),\dz\in(0,1)$, see Section \ref{notation}.

The corresponding global H\"older gradient estimate when $\partial\Omega,\varphi\in C^{2,\az}$ in the classical sense is given in the next theorem. 

\begin{thm}\label{s1: thm 2}
Let $u:\overline{\Omega}\to\mathbb{R}$ be a convex, continuous solution to \eqref{s1: eq 1}. Assume $\Omega$ is uniformly convex, $\partial\Omega,\varphi\in C^{2,\az}$ for some $\az\in(0,1)$. Then
$$u\in C^{1,\bz}(\overline{\Omega})$$
for some constant $\bz\in(0,1)$ depending only on $n,\lz,\Lz$ and $\az$.
\end{thm}

We will give an example to show that our results are optimal: if $\varphi$ is only $C^2$, the solution may fail to be globally $C^{1,\dz}$ for any $\dz\in(0,1)$.

Next we discuss case b) when the domain $\Omega$ is flat in a neighborhood of a boundary point. We have the following pointwise $C^{1,\az}$ estimate at a boundary point.

\begin{thm}\label{s1: thm 3}
Let $u:\overline{\Omega}\to\mathbb{R}$ be a convex, continuous solution to \eqref{s1: eq 1} with $\Omega=B^+_1$. Assume $\varphi\in C^{1,\az}(0)$ with $\az>\f{1}{5}$ and
$$\varphi(0)=0,\quad\quad\nabla_{x'}\varphi(0)=0,$$
and $\varphi$ separates quadratically on $\partial B^+_1$ in a neighborhood of $\{x_n=0\}$ from $0$. Then 
$$u\in C^{1,\az'}(0)$$
for some $\az'>0$ depending only on $n$ and $\az$. 
\end{thm}

The H\"older gradient estimate near the boundary in the flat case is as follows. 

We denote by $B'_R$ the ball in $\mathbb{R}^{n-1}$ centered at $0$ with radius $R>0$.

\begin{thm}\label{s1: thm 4}
Let $u:\overline{\Omega}\to\mathbb{R}$ be a convex, Lipschitz continuous solution to \eqref{s1: eq 1} with $\Omega=B^+_1$. Assume $\varphi|_{x_n=0}\in C^{1,\az}(B'_{3/4})$ with $\az>\max\{\f{1}{5},1-\f{2}{n}\}$, and for any $x'_0\in B'_{3/4}$, $\varphi|_{x_n=0}$ separates quadratically on $B'_1$ from its tangent plane at $x'_0$. Then 
$$u\in C^{1,\bz}(\overline{B^+_{1/2}})$$
for some small constant $\bz\in(0,1)$.
\end{thm}

In the particular case $\az=1$, Theorems \ref{s1: thm 3} and \ref{s1: thm 4} can be obtained from the work of the first author in \cite{S2} and \cite[Proposition 2.6]{S1}. The novelty here is that they hold when $\az < 1$.

The paper is organized as follows. In Section \ref{notation} we introduce some notation and give the quantitative versions of Theorems \ref{s1: thm 1}-\ref{s1: thm 4} (see Theorems \ref{thm 1}-\ref{thm 4} respectively). Section \ref{s3} is devoted to the proof of Theorem \ref{thm 1}. In Section \ref{s4}, we give the proof of Theorem \ref{thm 2}, and then present an example which shows that the assumptions in Theorem \ref{s1: thm 2} are sharp. In Sections \ref{s5} and \ref{s6}, we give the proofs of Theorems \ref{thm 3} and \ref{thm 4} respectively.

\section{Statement of main results}\label{notation}

We introduce some notation. We denote points in $\mathbb{R}^n$ as 
$$x=(x_1,\dots,x_n)=(x',x_n),\quad x'\in\mathbb{R}^{n-1}.$$
Let $u$ be a convex function defined on a convex set $\overline{\Omega}$, we 
denote by
$$l_{x_0}:=u(x_0)+\nabla u(x_0)\cdot(x-x_0)$$
a supporting hyperplane for the graph of $u$ at $x_0$ and $S_h(x_0)$ the section centered at $x_0$ and at height $h>0$,
$$S_{h}(x_0):=\{x\in\overline{\Omega}|\; u(x)<l_{x_0}(x)+h\}.$$
When $x_0\in\partial\Omega$, the term $\nabla u(x_0)$ is understood in the sense that
$$x_{n+1}=u(x_0)+\nabla u(x_0)\cdot(x-x_0)$$
is a supporting hyperplane for the graph of $u$ at $x_0$ but for any $\ez>0$,
$$x_{n+1}=u(x_0)+(\nabla u(x_0)+\ez\nu_{x_0})\cdot(x-x_0)$$
is not a supporting hyperplane, where $\nu_{x_0}$ denotes the unit inner normal to $\partial\Omega$ at $x_0$. We denote for simplicity $S_{h}=S_{h}(0)$, and sometimes when we specify the dependence on the function $u$ we use the notation $S_{h}(u)=S_{h}$. 

We state a variant of John's lemma \cite{Joh} (see also \cite{dG}), which is a classical result in convex geometry. 

\begin{lem}(See \cite{G}.)\label{Jo}
If $\Omega\subset\mathbb{R}^n$ is a bounded convex with nonempty interior and $E$ is the ellipsoid of minimum volume containing $\Omega$ centered at the center of mass of $\Omega$, then
$$\az_n E\subset\Omega\subset E,$$
where $\az_n=n^{-3/2}$ and $\az E$ denotes the $\az$-dilation of $E$ with respect to its center. 
\end{lem}

The following definition is introduced in \cite{S2}.

\begin{defn}\label{notation: def 1}
Let $k\ge 0$ be an integer and $0<\az\le 1$. We say that a function $u$ is pointwise $C^{k,\az}$ at $x_0$ and write
$$u\in C^{k,\az}(x_0)$$
if there exists a polynomial $P_{x_0}$ of degree $k$ such that
$$u(x)=P_{x_0}(x)+O(|x-x_0|^{k+\az}).$$
We say that $u\in C^k(x_0)$ if 
$$u(x)=P_{x_0}(x)+o(|x-x_0|^k).$$
\end{defn}

We now state the precise quantitative versions of Theorems \ref{s1: thm 1}-\ref{s1: thm 4} as follows.

\begin{thm}\label{thm 1}
Let $u:\overline{\Omega}\to\mathbb{R}$ be a convex, continuous solution to \eqref{s1: eq 1}. Assume $\Omega\subset\mathbb{R}^n_+, 0\in\partial\Omega$, $\Omega$ is uniformly convex at $0$, and on $\partial\Omega$ near $0$ we have
$$|x_n-q(x')|\le M|x'|^{2+\az},$$
$$|\varphi(x)-p(x')|\le M|x'|^{2+\az},$$
where $p(x'),q(x')$ are quadratic polynomials and 
$$M\ge\max\{\|\nabla p(0)\|,\,\|D^2_{x'}p\|,\,\|D^2_{x'}q\|\}.$$ 
Then 
$$u-u(0)-\nabla u(0)\cdot x\le C|x|^{1+\dz_0},$$
where $\dz_0>0$ depends only on $n$ and $\az$, the constant $C>0$ depends only on $n,\lz,\Lz,\az,M$, the uniform convexity of $\partial\Omega$ at $0$, and $\|\varphi\|_{L^{\infty}(\partial\Omega)}$. 
\end{thm}

This pointwise estimate combined with the interior estimates of Caffarelli from \cite{C3} implies the global $C^{1,\alpha}$ estimate for solutions to \eqref{s1: eq 1} in the case that the domain is uniformly convex.

\begin{thm}\label{thm 2}
Let $u:\overline{\Omega}\to\mathbb{R}$ be a convex, continuous solution to \eqref{s1: eq 1}. Assume $\Omega$ is uniformly convex, $\partial\Omega,\varphi\in C^{2,\az}$ for some $\az\in(0,1)$. Then
$$[\nabla u]_{C^{\bz}(\overline{\Omega})}\le C,$$
where $\bz\in(0,1)$ depends only on $n,\lz,\Lz$ and $\az$, the constant $C>0$ depends only on $n,\lz,\Lz,\az$, $\mathrm{diam}(\Omega), \|\partial\Omega,\varphi\|_{C^{2,\az}}$ and the uniform convexity of $\Omega$.
\end{thm}

In the case that the domain is flat at a boundary point, the quantitative pointwise $C^{1,\az}$ estimate is as follows.

\begin{thm}\label{thm 3}
Let $u:\overline{\Omega}\to\mathbb{R}$ be a convex, continuous solution to \eqref{s1: eq 1} with $\Omega=B^+_1$ and
$$u(0)=0,\quad\quad\nabla u(0)=0.$$
Assume $\az>\f{1}{5}$, and
$$\varphi|_{x_n=0}\le\mu^{-1}|x'|^{1+\az}\quad\mathrm{in}\;B'_{1/2}$$ 
and 
\begin{equation}\label{s2: eq 3.4}
\varphi\ge\mu|x|^{2}\quad\mathrm{on}\;\partial B^+_1\cap\{x_n\le\rz\}.
\end{equation}
for some $\mu,\rz>0$. Then for any $x\in B^+_1$ with $u(x)\le c$, we have
\begin{equation*}
u(x)\le C|x|^{1+\az'},
\end{equation*}
where $\az'>0$ depends only on $n$ and $\az$, the constants $c, C$ depend only on $n,\lz,\Lz,\mu,\az$ and $\rz$.
\end{thm}

Using Theorem \ref{thm 3} and similar techniques as in the uniformly convex case, we can obtain the $C^{1,\az}$ estimate near the flat boundary.

\begin{thm}\label{thm 4}
Let $u:\overline{\Omega}\to\mathbb{R}$ be a convex, Lipschitz continuous solution to \eqref{s1: eq 1} with $\Omega=B^+_1$. Assume $\varphi|_{x_n=0}\in C^{1,\az}(B'_{3/4})$ with $\az>\max\{\f{1}{5},1-\f{2}{n}\}$, and for any $x_0=(x'_0,0)$ with $x'_0\in B'_{3/4}$ and $x\in\partial B^+_1\cap\{x_n=0\}$, 
$$\varphi(x)-\varphi(x_0)-\nabla_{x'}\varphi(x_0)\cdot(x'-x'_0)\ge\mu|x'-x'_0|^2.$$
Then 
$$[\nabla u]_{C^{\bz}(\overline{B^+_{1/2}})}\le C,$$
where $\bz\in(0,1)$ depends only on $n,\lz,\Lz$,$\az$, and the constant $C>0$ depends on $n,\lz,\Lz,\az,\mu$, $\|\varphi|_{x_n=0}\|_{C^{1,\az}(B'_{3/4})}$ and $\|u\|_{C^{0,1}}$. 
\end{thm}

In the proofs below we denote by $c, C, c', C', c_i, C_i (i=0,1,2,\dots,)$ constants depending only on the data $n,\lz,\Lz,\az,\mathrm{diam}(\Omega),M$, $\|\varphi\|_{L^{\infty}(\partial\Omega)}$, the uniform convexity of $\partial\Omega$ etc. Their values may change from line to line whenever there is no possibility of confusion.
 For $A,B\in\mathbb{R}$, we write $A\sim B$ if 
$$c\le\f{A}{B}\le C$$
for some universal constants $c, C$.

\section{Proof of Theorem \ref{thm 1}}\label{s3}

Let $\varphi(x)=\tilde{\varphi}(x')$ and $\tilde{u}=u-l_0$, where we recall from Section \ref{notation} that 
$$l_0(x)=u(0)+\nabla u(0)\cdot x.$$
Then (after performing a rotation in the $x'$ subspace) on $\partial\Omega$ we have
$$\tilde{u}=\tilde{\varphi}-\tilde{\varphi}(0)-\nabla_{x'}\tilde{\varphi}(0)\cdot x'-u_n(0)x_n=\sum_1^{n-1}a_i^2 x_i^2+O(|x'|^{2+\az})$$
for some constants $a_i\ge 0, i=1,\dots,n-1$.

Let $0<\az'<\az$ be a constant to be chosen below. We will prove that 
\begin{equation}\label{s5: thm 20}
S_h\supset\overline{\Omega}\cap B_{ch^{\f{1}{1+\dz_0}}}\quad\forall\;h>0,
\end{equation}
where $\dz_0>0$ is a constant depending only on $n$ and $\az$.

First, we use a lower barrier of the type
$$\tilde{\varphi}(0)+\nabla_{x'}\tilde{\varphi}(0)\cdot x'+\Lz|x|^2-Cx_n$$
and obtain that $u_n(0)$ is bounded. Hence $\tilde{u}$ is bounded above and therefore we can assume that $h$ in \eqref{s5: thm 20} is sufficiently small.

We only need to consider the following cases: $\min_ia_i^2\le h^{\f{\az'}{2+\az'}}$ and $\min_ia_i^2\ge h^{\f{\az'}{2+\az'}}$.

\noindent$\mathbf{Case\;1:}$ $\min_ia_i^2\le h^{\f{\az'}{2+\az'}}$. 

If $a_1^2\le h^{\f{\az'}{2+\az'}}$, then by the uniform convexity of $\partial\Omega$ at $0$, we have on $\partial\Omega\cap\{x_i=0,i=2,\dots,n-1\}$
$$\tilde{u}\le h^{\f{\az'}{2+\az'}}x_1^2+O(|x_1|^{2+\az})\le C\l[h^{\f{\az'}{2+\az'}}x_n+x_n^{1+\f{\az'}{2}}\r],$$
this together with the convexity of $u$ implies that 
$$\tilde{u}(te_n)\le C\l[h^{\f{\az'}{2+\az'}}t+t^{1+\f{\az'}{2}}\r].$$
It follows that 
$$\{te_n: 0\le t\le c_0h^{\f{2}{2+\az'}}\}\subset S_h$$
for some small constant $c_0>0$. 

The domain of definition of $\partial S_h\cap\partial\Omega$ contains a ball in $\mathbb{R}^{n-1}$ of radius $ch^{\f{1}{2}}$, and by the uniform convexity of $\partial\Omega$ at $0$, we have 
$$x_n\ge c_1h\quad\quad\forall\,x\in\partial S_h\cap\partial\Omega\cap\{|x'|=ch^{\f{1}{2}}\}.$$
and therefore
\begin{equation*}
S_h\supset\overline{\Omega}\cap\{x_n\le c_1h\}.
\end{equation*}
Then the convex set generated by $\Omega\cap\{x_n=c_1h\}$ and the point $(0,c_0h^{\f{2}{2+\az'}})$ is contained in $S_h$. Since this convex set contains a half-ball centered at $(0, c_1h)$ of radius $ch^{\f{2}{2+\az'}}$. We obtain that
$$S_h\supset\overline{\Omega}\cap B_{ch^{\f{2}{2+\az'}}}.$$

\noindent$\mathbf{Case\;2:}$ $\min_ia_i^2\ge h^{\f{\az'}{2+\az'}}$. Then on $\partial\Omega$ near $0$ we have 
\begin{equation}\label{s5: thm 21}
\tilde{u}\ge\f{1}{2n}h^{\f{\az'}{2+\az'}}|x'|^2\quad\quad\forall x\in\{|x'|\le c'h^{\f{\az'}{\az(2+\az')}}\},
\end{equation}
for some $c'$ small.

Let $x^*_h$ be the center of mass of $S_h$ and denote $d_h:=x^*_h\cdot e_n$. We claim that
\begin{equation}\label{s5: thm 22}
d_h\ge c_2h^{\f{n}{n+1}}
\end{equation}
for some constant $c_2$ small. Otherwise, by the uniformly convexity of $\partial\Omega$ at $0$ and Lemma \ref{Jo}, we have
$$S_h\subset\{0\le x_n\le C(n)c_2h^{\f{n}{n+1}}\le h^{\f{n}{n+1}}\}\cap\{|x'|\le C_1h^{\f{n}{2(n+1)}}\}.$$
Let $C_2$ be a large constant to be chosen and define
$$w:=\ez x_n+\l[\f{1}{2}\l(\f{|x'|}{C_1h^{\f{n}{2(n+1)}}}\r)^2+C_2\l(\f{x_n}{h^{\f{n}{n+1}}}\r)^2\r]\cdot h.$$
Since 
$$\partial S_h\cap\partial\Omega\subset\{|x'|\le C_1h^{\f{n}{2(n+1)}}\}\subset\{|x'|\le c'h^{\f{\az'}{\az(2+\az')}}\}
$$
if we choose $\az'$ small, then on $\partial S_h\cap\partial\Omega$ we have
$$w\le\l[C\ez+\f{1}{2C_1^2}h^{\f{1}{n+1}}+C_2Ch^{\f{1}{n+1}}\r]|x'|^2\le\tilde{u},$$
where we choose $\ez$ and $\az'$ small such that
$$C\ez\le\f{h^{\f{\az'}{2+\az'}}}{6n}\quad\mathrm{and}\quad C_2Ch^{\f{1}{n+1}}\le\f{h^{\f{\az'}{2+\az'}}}{6n}.$$
In $S_h$ we have
$$w\le\ez+\l[\f{1}{2}+C_2C(n)c_2\r]h\le h.$$
Moreover,
$$\mathrm{det}\;D^2 w>\Lz$$
by choosing $C_2$ large.

In conclusion, $w\le\tilde{u}$ in $S_h$, which together with the convexity of $u$ implies that $\tilde{u}\ge\ez x_n$ in $\Omega$. This is a contradiction. Thus \eqref{s5: thm 22} holds.

The uniform convexity of $\Omega$ at $0$ and \eqref{s5: thm 22} imply that 
$$\{te_n: 0\le t\le c_2h^{\f{n}{2(n+1)}+\f{1}{2}}\}\subset S_h$$
for some small constant $c_2>0$. Similar to $\mathbf{Case\;1}$ we have
$$S_h\supset\overline{\Omega}\cap B_{ch^{\f{n}{2(n+1)}+\f{1}{2}}}.$$

Combining $\mathbf{Cases\;1}$-$\mathbf{2}$, we obtain \eqref{s5: thm 20}.

\section{Proof of Theorem \ref{thm 2}}\label{s4}

Using the uniform pointwise estimate of Theorem \ref{thm 1} we obtain by standard arguments the H\"older continuity of $\nabla u$ on $\partial \Omega$.

\begin{lem}\label{lem 4.1}
Under the assumptions of Theorem \ref{thm 2}, we have
$$[\nabla u]_{C^{\dz_0}(\partial\Omega)}\le C,$$
where $\dz_0$ is the constant in Theorem \ref{thm 1}.

\end{lem}

It remains to show that $\nabla u$ is uniformly H\"older continuous also at interior points of $\Omega$.
Assume $0\in\partial\Omega$ and $\Omega\subset\mathbb{R}^n_+$. We divide the proof of Theorem \ref{thm 2} into three steps.\\

\noindent$\mathbf{Step\;1.}$ Let $y\in\Omega$ and consider the maximal interior section $S_{\bar{h}}(y)$ centered at $y$, that is,
$$\bar{h}=\max\{h|\;S_h(y)\subset\Omega\}.$$
Assume $0\in\partial S_{\bar{h}}(y)\cap\partial\Omega$. We prove that
\begin{equation}\label{thm 21}
S_{\bar{h}}(y)\subset B_{\bar{h}^{\ez_0}}
\end{equation}
for any $\bar{h}>0$ small, where $\ez_0$ is a small constant.

For any $h>0$ small, let $x^*_h$ be the center of mass of $S_h$ and $d_h:=x^*_h\cdot e_n$. We claim that
\begin{equation}\label{thm 20}
d_h\le C_0h^{\f{1}{2}}
\end{equation}
for some large constant $C_0>0$.

Indeed, if \eqref{thm 20} does not hold, then as in the proof of Theorem \ref{thm 1} we have
$$\overline{S_h}\supset\Omega\cap\{x_n\le c_1h\}$$
for some $c_1$ small, and therefore $\overline{S_h}$ contains the convex set generated by $\Omega\cap\{x_n=c_1h\}$ and the point $x^*_h$. We also have
$$|x'|\ge ch^{\f{1}{2}},\quad\forall\,x\in\partial\Omega\cap\{x_n=c_1h\}.$$ 
It follows that
$$|S_h|\ge c(n)(ch^{\f{1}{2}})^{n-1}\f{C_0h^{\f{1}{2}}}{2}.$$
On the other hand, by Lemma \ref{Jo}, $S_h$ is equivalent to an ellipsoid $E$ centered at $x^*_h$, i.e.,
$$E\subset S_h\subset C(n)E,$$
where the dilation is with respect to $x^*_h$. Let $P$ be the quadratic polynomial that solves
$$\mathrm{det}\;D^2 P=\lz\quad\mathrm{in}\;E,\quad\quad P=h\quad\mathrm{on}\;\partial E.$$
Then
$$P\ge\tilde{u}:=u-l_0\ge 0\quad\mathrm{in}\;E.$$
It follows 
$$h^n\ge|h-\min_{E}P|^n\ge c(n,\lz)|E|^2\ge c|S_h|^2.$$ 
We reach a contradiction if we choose $C_0$ sufficiently large. Hence \eqref{thm 20} holds, which gives
\begin{eqnarray}\label{thm 200}
S_h\subset\{0\le x_n\le C'h^{\f{1}{2}}\}\quad\quad\forall\,h>0.
\end{eqnarray}
This together with the uniform convexity of $\Omega$ gives
\begin{eqnarray}\label{thm 20'}
S_h\subset B_{Ch^{\f{1}{4}}}\quad\quad\forall\,h>0.
\end{eqnarray}

Since $\partial\Omega,\varphi\in C^2$, we find that
$$S_{\bar{h}}(y)=\{x\in\Omega|\;v(x)<0\},$$
where 
$$v(x)=(u-l_0)(x)-[u_n(y)-u_n(0)]x_n.$$
Choose $h>0$ such that 
$$u_n(y)-u_n(0)=(2C')^{-1}h^{\f{1}{2}},$$
where $C'$ is the constant in \eqref{thm 200}. Then we have by \eqref{thm 200} 
$$S_{\bar{h}}(y)\cap\{x_n\le 2C'h^{\f{1}{2}}\}\subset S_h\subset\{x_n\le C'h^{\f{1}{2}}\},$$
which implies that
$$S_{\bar{h}}(y)\subset\{x_n\le 2C'h^{\f{1}{2}}\}.$$
Using \eqref{thm 20'} we obtain that
\begin{equation}\label{thm 21'}
S_{\bar{h}}(y)\subset S_h\subset B_{Ch^{\f{1}{4}}}.
\end{equation}
On the other hand, let $\tz>0$ be a small constant to be chosen below and denote 
$$\Lz_0:=\f{1+\dz_0}{2\dz_0},$$
where $\dz_0$ is the constant in Theorem \ref{thm 1}. By Theorem \ref{thm 1}, we can choose a point $z=te_n\in\partial S_{\tz h^{\Lz_0}}$ with $t\ge c(\tz h^{\Lz_0})^{\f{1}{1+\dz_0}}$. It follows that
\begin{eqnarray*}
v(z)&\le&\tz h^{\Lz_0}-(2C')^{-1}h^{\f{1}{2}}c(\tz h^{\Lz_0})^{\f{1}{1+\dz_0}}\\
&=&\tz^{\f{1}{1+\dz_0}}h^{\Lz_0}[\tz^{\f{\dz_0}{1+\dz_0}}-(2C')^{-1}c]\\
&\le&-\tz^{\f{1}{1+\dz_0}}h^{\Lz_0}(4C')^{-1}c<0
\end{eqnarray*}
if $\tz>0$ is sufficiently small. This implies that
\begin{equation}\label{thm 21''}
\bar{h}\ge c|S_{\bar{h}}(y)|^{\f{2}{n}}\ge c|\min_{S_{\bar{h}}(y)}v|\ge ch^{\Lz_0}.
\end{equation}
This together with \eqref{thm 21'} gives \eqref{thm 21} with $0<\ez_0<\f{1}{4\Lz_0}$.\\

\noindent$\mathbf{Step\;2.}$ Let $S_{\bar{h}}(y)$ be a maximal interior section tangent to $\partial\Omega$ at $0$ as in $\mathbf{Step\;1}$. We prove that
\begin{equation}\label{thm 21**}
B_{\bar{h}^{1-\ez_1}}(y)\subset S_{\bar{h}}(y)
\end{equation}
for any $\bar{h}>0$ small, where $\ez_1>0$ is a small constant.

By Lemma \ref{Jo}, $S_{\bar{h}}(y)$ is equivalent to an ellipsoid $E$ centered at $y$, i.e., 
$$E\subset S_{\bar{h}}(y)\subset C_1E,$$
where $C_1$ is a constant depending only on $n,\lz,\Lz$, and
\begin{eqnarray*}
E=y+U^t\mathrm{diag}(\mu_1,\dots,\mu_n)B_1,
\end{eqnarray*}
where $0<\mu_1\le\mu_2\le\dots\le\mu_n$ and $U$ is an orthogonal matrix.  We only need to prove that
\begin{equation}\label{thm 22}
\mu_1\ge\bar{h}^{1-\ez_1}.
\end{equation}
Assume by contradiction that 
$$\mu_1<\bar{h}^{1-\ez_1}.$$
Let $\nu=U^te_1$ be a unit vector which is parallel to the shortest axis of $E$. Then for any $x\in C_1E$, we have
\begin{eqnarray*}
|(x-y)\cdot\nu|&=&|e_1^t\mathrm{diag}(\mu_1,\dots,\mu_n)\mathrm{diag}(\mu_1^{-1},\dots,\mu_n^{-1})U(x-y)|\\
&\le&C_1|e_1^t\mathrm{diag}(\mu_1,\dots,\mu_n)|\\
&=&C_1\mu_1\le C_1\bar{h}^{1-\ez_1}.
\end{eqnarray*}
Define $w^+=v-\f{\bar{h}^{\ez_1}}{2C_1}(x-y)\cdot\nu$ and $a^+=\min_{\overline{S_{\bar{h}}(y)}}w^+=w^+(x_0)$. Since
$$w^+(y)=v(y)=-\bar{h},\quad\quad w^+=-\f{\bar{h}^{\ez_1}}{2C_1}(x-y)\cdot\nu\ge-\f{\bar{h}}{2}\quad\mathrm{on}\;\partial S_{\bar{h}}(y),$$
we find that $x_0\in S_{\bar{h}}(y)$ and
$$v\ge\f{\bar{h}^{\ez_1}}{2C_1}(x-y)\cdot\nu+w^+(x_0)\quad\mathrm{in}\;S_{\bar{h}}(y).$$
It follows from the convexity of $v$ that
$$v\ge\f{\bar{h}^{\ez_1}}{2C_1}(x-y)\cdot\nu+w^+(x_0)\ge\f{\bar{h}^{\ez_1}}{2C_1}(x-y)\cdot\nu-2\bar{h}\quad\mathrm{in}\;\Omega.$$
Similarly we have
$$v\ge\f{\bar{h}^{\ez_1}}{2C_1}(x-y)\cdot(-\nu)-2\bar{h}\quad\mathrm{in}\;\Omega.$$
The last two estimates imply 
\begin{equation}\label{thm 24'}
v\ge c\bar{h}^{\ez_1}|(x-y)\cdot\nu|-2\bar{h}\quad\mathrm{in}\;\Omega.
\end{equation}

Recall that
$$v=u-l_y-\bar{h}=u-l_0-[u_n(y)-u_n(0)] x_n$$
satisfies
\begin{equation}\label{thm 23}
v<0\quad\mathrm{in}\;S_{\bar{h}}(y),\quad\quad v\ge 0\quad\mathrm{in}\;\Omega\setminus S_{\bar{h}}(y).
\end{equation}
Since $v\ge 0$ on $\partial\Omega$ and $\partial\Omega,\varphi\in C^{2,\az}$, we have (after performing a rotation in the $x'$ subspace)
\begin{equation}\label{thm 24''}
v|_{\partial\Omega}=\sum_{1}^{n-1}\lz_i^2x_i^2+O(|x'|^{2+\az})
\end{equation}
for some bounded constants $\lz_i,i=0,\dots,n-1$.

We only need to consider the cases: $|\nu'|\ge\bar{h}^{\f{\ez_0}{4}}$ and $|\nu'|\le\bar{h}^{\f{\ez_0}{4}}$, where $\ez_0$ is the constant in \eqref{thm 21}. 

\noindent$\mathbf{Case\;1.}$ $|\nu'|\ge\bar{h}^{\f{\ez_0}{4}}$. We choose $x\in\partial\Omega$ with $x'=\bar{h}^{\f{\ez_0}{2}}\f{\nu'}{|\nu'|}$. Then we have
$$x_n\le C|x'|^2\le C\bar{h}^{\f{\ez_0}{2}}|x'|,$$
and
$$|y|\le\bar{h}^{\ez_0}\le\bar{h}^{\f{\ez_0}{2}}|x'|.$$
It follows that 
\begin{eqnarray*}
\bar{h}^{\ez_1}(x-y)\cdot\nu&\ge&\bar{h}^{\ez_1}\l(x'\cdot\nu'-x_n-|y|\r)\\
&\ge&\bar{h}^{\ez_1}\l(|\nu'||x'|-C\bar{h}^{\f{\ez_0}{2}}|x'|\r)\\
&\ge&\bar{h}^{\ez_1+\f{\ez_0}{4}}(1-C\bar{h}^{\f{\ez_0}{4}})|x'|\\
&\ge&\f{\bar{h}^{\ez_1+\f{\ez_0}{4}}}{2}|x'|
\end{eqnarray*}
if $\ez_1,\ez_0>0$ are sufficiently small and $\bar{h}$ is small. It follows from \eqref{thm 24'} and \eqref{thm 24''} that
\begin{eqnarray*}
c\bar{h}^{\ez_1}\bar{h}^{\f{\ez_0}{4}}|x'|\le v(x)\le C|x'|^2=C\bar{h}^{\f{\ez_0}{2}}|x'|.
\end{eqnarray*}
Choose $0<\ez_1<\f{\ez_0}{4}$ and then we reach a contradiction. Thus \eqref{thm 22} holds.

\noindent$\mathbf{Case\;2.}$ $|\nu'|\le\bar{h}^{\f{\ez_0}{4}}$. Then we have
$$|\nu\cdot e_n|\ge 1-\bar{h}^{\f{\ez_0}{4}}>\f{1}{2}$$
if $\bar{h}$ is small.

For any $x\in\Omega$ near $0$ with $x_n^{\f{1}{2}}\ge\bar{h}^{\f{\ez_0}{8}}$, we have
\begin{eqnarray*}
\bar{h}^{\ez_1}|(x-y)\cdot\nu|&\ge&\bar{h}^{\ez_1}\l(x_n|\nu\cdot e_n|-|x'||\nu'|-|y|\r)\\
&\ge&\bar{h}^{\ez_1}\l(\f{1}{2}x_n-C\bar{h}^{\f{\ez_0}{4}}x_n^{\f{1}{2}}-\bar{h}^{\ez_0}\r)\\
&\ge&\f{\bar{h}^{\ez_1}}{4}x_n
\end{eqnarray*}
if $\ez_1,\ez_0>0$ are sufficiently small and $\bar{h}$ is small. Hence by \eqref{thm 24'},
\begin{equation}\label{thm 27}
v\ge c\bar{h}^{\ez_1}x_n,\quad\quad\mathrm{in}\;\Omega\cap\{\bar{h}^{\f{\ez_0}{8}}\le x_n^{\f{1}{2}}\le c\}.
\end{equation}

Let $0<\dz<\f{\ez_0}{8}$ be a small constant to be chosen below.

\noindent$\mathbf{Case\;2.1.}$ If one of $\lz_i^2, i=1,\dots, n-1$, say $\lz_1^2$, satisfies $\lz_1^2\le\bar{h}^{\dz\az}$, then we choose $x=(x_1,0,\dots,0,x_n)\in\partial\Omega$ with $x_n^{\f{1}{2}}=\bar{h}^{\dz}$. We have by \eqref{thm 24''} and \eqref{thm 27}
$$c\bar{h}^{\ez_1}x_n\le v(x)\le\lz_1^2x_1^2+C|x_1|^{2+\az}\le C\bar{h}^{\dz\az}x_n.$$
Choose $\ez_1<\dz\az$ and we reach a contradiction.

\noindent$\mathbf{Case\;2.2.}$ $\min_{1\le i\le n-1}\lz_i^2\ge\bar{h}^{\dz\az}$. Then we have
\begin{equation}\label{thm 27'}
v\ge\f{1}{2n}\bar{h}^{\dz\az}|x'|^2,\quad\quad\mathrm{on}\;\partial\Omega\cap\{x_n^{\f{1}{2}}\le c\bar{h}^{\dz}\}.
\end{equation}
Define
$$w=\sz x_n+\f{1}{4n}\bar{h}^{\dz\az}|x'|^2+\f{C_*}{\bar{h}^{\dz\az(n-1)}}x_n^2.$$
Then $w$ is a lower barrier for $v$ in $\Omega\cap\{x_n^{\f{1}{2}}\le\bar{h}^{\f{\ez_0}{8}}\}$ if $\sz,\dz>0$ are sufficiently small and $C_*$ is large. 

Indeed, on $\partial\Omega\cap\{x_n^{\f{1}{2}}\le\bar{h}^{\f{\ez_0}{8}}\}$ we have 
$$w\le C\sz|x'|^2+\f{1}{4n}\bar{h}^{\dz\az}|x'|^2+C\bar{h}^{\f{\ez_0}{4}-\dz\az(n-1)}|x'|^2\le\f{1}{2n}\bar{h}^{\dz\az}|x'|^2$$
if $\sz$ is small and $\dz\az n<\f{\ez_0}{4}$.

On $\Omega\cap\{x_n^{\f{1}{2}}=\bar{h}^{\f{\ez_0}{8}}\}$ we have 
$$w\le\sz x_n+C\bar{h}^{\dz\az}x_n+C\bar{h}^{\f{\ez_0}{4}-\dz\az(n-1)}x_n\le c\bar{h}^{\ez_1}x_n$$
if we choose $\ez_1<\min\l\{\dz\az,\f{\ez_0}{4}-\dz\az(n-1)\r\}$ and $\sz$ small.

Hence by \eqref{thm 27'} and \eqref{thm 27} we obtain that $v\ge w\ge\sz x_n$ in $\Omega\cap\{x_n^{\f{1}{2}}\le\bar{h}^{\f{\ez_0}{8}}\}$. Since $y\in\Omega\cap B_{\bar{h}^{\ez_0}}\subset\Omega\cap\{x_n^{\f{1}{2}}\le\bar{h}^{\f{\ez_0}{8}}\}$, we reach a contradiction since $v(y)=-\bar{h}<0$.

Combining $\mathbf{Case\;2.1}$ and $\mathbf{Case\;2.2}$, we prove \eqref{thm 22} in $\mathbf{Case\;2}$. Hence \eqref{thm 21**} holds.\\

\noindent$\mathbf{Step\;3.}$ We show that
\begin{equation}\label{thm 20*}
|\nabla u(x)-\nabla u(y)|\le C|x-y|^{\bz}\quad\quad\forall\,x,y\in\overline{\Omega},
\end{equation}
where $\bz\in(0,1)$ is a constant depending only on $n,\lz,\Lz$ and $\az$. 

This follows from $\mathbf{Steps\;1}$-$\mathbf{2}$ and similar arguments as in \cite{S1}. For completeness, we include the proof.

We first note that in $\mathbf{Steps\;1}$-$\mathbf{2}$, if $\bar{h}\ge c$ for some small constant $c$, the estimates \eqref{thm 21} and \eqref{thm 21**} obviously hold since $\bar{h}\sim|S_{\bar{h}}(y)|^{\f{2}{n}}$ is bounded above. (We only need to replace $\bar{h}^{1-\ez_1}$ by $C^{-1}\bar{h}^{1-\ez_1}$ and $\bar{h}^{\ez_0}$ by $C\bar{h}^{\ez_0}$ for some large constant $C$.)

Let $y\in\Omega$ and assume the maximal interior section $S_{\bar{h}}(y)$ is tangent to $\partial\Omega$ at $0\in\partial\Omega$. Let $Tx=Ax+b$ be an affine transformation that normalizes $S_{\bar{h}}(y)$, i.e.,
$$B_{\az_n}\subset TS_{\bar{h}}(y)\subset B_1.$$
By \eqref{thm 21**} we have
\begin{equation}\label{thm 21*}
\|T\|\le C\bar{h}^{-(1-\ez_1)},\quad\quad |\mathrm{det}\;T|^{\f{2}{n}}\sim |S_{\bar{h}}(y)|^{-\f{2}{n}}\sim\bar{h}^{-1}.
\end{equation}
For $\tilde{x}\in TS_{\bar{h}}(y)$, define
$$\tilde{u}(\tilde{x})=|\mathrm{det}\;T|^{\f{2}{n}}[u-l_y-\bar{h}](T^{-1}\tilde{x}),$$
where we recall that
$$l_y(x)=u(y)+\nabla u(y)\cdot(x-y).$$
The interior $C^{1,\gz}$ estimate for solutions of the Monge-Amp$\grave{e}$re equation (see \cite{G}) gives
$$|\nabla\tilde{u}(\tilde{x}_1)-\nabla\tilde{u}(\tilde{x}_2)|\le C|\tilde{x}_1-\tilde{x}_2|^{\gz}\quad\forall\,\tilde{x}_1,\tilde{x}_2\in TS_{\f{\bar{h}}{2}}(y)$$
for some $\gz\in(0,1)$ depending only on $n,\lz,\Lz$. Rescaling back and using
$$\nabla\tilde{u}(\tilde{x}_1)-\nabla\tilde{u}(\tilde{x}_2)=|\mathrm{det}\;T|^{\f{2}{n}}[\nabla u(T^{-1}\tilde{x}_1)-\nabla u(T^{-1}\tilde{x}_2)]A^{-1},$$
we find from \eqref{thm 21*}
\begin{eqnarray}\label{thm 22*}
|\nabla u(x_1)-\nabla u(x_2)|&\le&C|\mathrm{det}\;T|^{-\f{2}{n}}\|A\|^{1+\gz}|x_1-x_2|^{\gz}\nonumber\\
&\le&C\bar{h}^{1-(1-\ez_1)(1+\gz)}|x_1-x_2|\nonumber\\
&\le&C|x_1-x_2|^{\gz}\quad\quad\forall\,x_1,x_2\in S_{\f{\bar{h}}{2}}(y).
\end{eqnarray}
if we choose $\gz>0$ sufficiently small.

The convexity of $u$ implies that $S_{\f{\bar{h}}{2}}(y)\supset\f{1}{2}S_{\bar{h}}(y)$, where the rescaling is with respect to $y$. Hence by \eqref{thm 21**},
$$S_{\f{\bar{h}}{2}}(y)\supset B_{c\bar{h}^{1-\ez_1}}(y).$$

Also, by the proof of \eqref{thm 21} (see \eqref{thm 21''}) we have
\begin{equation}\label{thm 24*}
|\nabla u(y)-\nabla u(0)|\le C\bar{h}^{2\ez_0}.
\end{equation}

For any $x,y\in\Omega$, assume the maximal interior sections $S_{\bar{h}_x}(x), S_{\bar{h}_y}(y)$ are tangent to $\partial\Omega$ at some points $\bar{x}, \bar{y}\in\partial\Omega$ respectively. Assume without loss of generality that $\bar{h}_y\ge\bar{h}_x$.

\noindent$\mathbf{Case\;1.}$ $x\in S_{\f{\bar{h}_y}{2}}(y)$. Then by \eqref{thm 22*} we have
\begin{equation*}
|\nabla u(x)-\nabla u(y)|\le C|x-y|^{\gz}.
\end{equation*}

\noindent$\mathbf{Case\;2.}$ $x\notin S_{\f{\bar{h}_y}{2}}(y)$. Then we have
$$|x-y|\ge c\bar{h}_y^{1-\ez_1}\ge c\bar{h}_x^{1-\ez_1}.$$
and therefore by \eqref{thm 21}
$$|\bar{x}-\bar{y}|\le |\bar{x}-x|+|x-y|+|y-\bar{y}|\le C[\bar{h}_x^{\ez_0}+|x-y|+\bar{h}_y^{\ez_0}]\le C|x-y|^{\f{\ez_0}{1-\ez_1}}.$$
Hence by \eqref{thm 24*} and Lemma \ref{lem 4.1},
\begin{eqnarray*}
|\nabla u(x)-\nabla u(y)|&\le&|\nabla u(x)-\nabla u(\bar{x})|+|\nabla u(\bar{x})-\nabla u(\bar{y})|+|\nabla u(\bar{y})-\nabla u(y)|\\
&\le&C[\bar{h}_x^{2\ez_0}+|\bar{x}-\bar{y}|^{\dz_0}+\bar{h}_y^{2\ez_0}]\\
&\le&C|x-y|^{\f{\ez_0\dz_0}{1-\ez_1}}.
\end{eqnarray*}
Hence we obtain \eqref{thm 20*}. The proof of Theorem \ref{thm 2} is complete.\\

\qed

We conclude this section with an example which shows that if the boundary data $\varphi$ is only $C^2$ in Theorem \ref{thm 2}, then $u$ may fail to be of class $C^{1,\dz}(\overline{\Omega})$ for any $\dz\in(0,1)$.

\bigskip

\noindent $\mathbf{Example.}$ 
Let $\Omega=B_{\rz}(\rz e_n)$, where $\rz$ is small depending only on $n$ to be chosen below. Let $u$ solves
$$\l\{
\begin{array}{rclcl}
\mathrm{det}\;D^2u&=&1&&\mathrm{in}\;\Omega,\\
u&=&\f{x_n}{-\log x_n}&&\mathrm{on}\;\partial\Omega,
\end{array}
\r.$$
where we define $u(0)=\lim_{|x|\to 0}\f{x_n}{-\log x_n}=0$.

In a neighborhood of $0$, the boundary data $\varphi=u|_{\partial\Omega}$ can be written as 
$$\varphi(x')=\f{\rz-\sqrt{\rz^2-|x'|^2}}{-\log(\rz-\sqrt{\rz^2-|x'|^2})}.$$
By straightforward computation we find that $\varphi\in C^2(\partial\Omega)$.

Next we show that $u_n(0)\le 0$. 

Indeed, for any $0<t<\rz$, we choose $y=(y',t)\in\partial\Omega$. Then the convexity of $u$ gives
$$\f{u(te_n)}{t}\le\f{1}{2}\l[\f{u(y',t)}{t}+\f{u(-y',t)}{t}\r]=\f{1}{-\log t},$$
which implies $u_n(0)\le 0$.

Now we construct a lower barrier for $u$ in $\Omega$.

Let
$$w:=\f{1}{2}\f{x_n}{-\log x_n}+|x'|^2x_n^{\f{1}{n}}.$$
Then we can compute 
$$\mathrm{det}\;D^2 w=2^{n-2}\l[\f{1}{x_n^{\f{1}{n}}(\log x_n)^2}\l(1+\f{2}{-\log x_n}\r)-\f{2(n+1)}{n^2}|x'|^2x_n^{-1}\r]\quad\mathrm{in}\;\Omega.$$
Since $|x'|^2x_n^{-1}$ and $x_n$ are bounded by $2\rz$, we can choose $\rz>0$ small depending only on $n$ such that 
$$\mathrm{det}\;D^2w>1\quad\mathrm{in}\;\Omega$$
and
$$w\le\f{1}{2}\f{x_n}{-\log x_n}+2\rz x_n^{1+\f{1}{n}}\le\f{x_n}{-\log x_n}\quad\mathrm{on}\;\partial\Omega.$$
Therefore $u\ge w\ge 0$ in $\Omega$, which implies $u_n(0)\ge 0$. Hence, 
$$u_n(0)=0.$$
It follows that
$$u(0,x_n)-u_n(0)x_n=u(0,x_n)\ge w(0,x_n)=\f{1}{2}\f{x_n}{-\log x_n}.$$
This implies that $u\notin C^{1,\dz}(0)$ for any $\dz\in(0,1)$.

\section{Proof of Theorem \ref{thm 3}}\label{s5}

We divide the proof into two steps.

\

\noindent$\mathbf{Step\;1:}$ We construct an explicit barrier for $u$. Let
$$\bar{w}(r,y):=r^{2}(1-t^{\ez})^+,\quad t=yr^{-\sz}\ge 0,$$
for some $\sz>0,\ez\in(0,1)$ to be chosen below. Then the function
$$w_1(x',x_n):=c'\bar{w}(|x'|, C'x_n),$$
is a lower barrier for $u$ provided that $c'$ (small), $C'$ (large) are  appropriate constants.

Since
$$\f{dt}{dr}=-\sz r^{-1}t,\quad\quad\f{dt}{dy}=r^{-\sz},$$
we compute in the set $B^+_1\cap\{\bar{w}>0\}$ (i.e. $t\in(0,1)$):
\begin{eqnarray*}
\bar{w}_r&=&2r(1-t^{\ez})+\sz\ez rt^{\ez},\\
\bar{w}_{rr}&=&2(1-t^{\ez})+(3-\sz\ez)\sz\ez t^{\ez},\\
\bar{w}_{ry}&=&(\sz\ez-2)\ez t^{\ez-1}r^{1-\sz},\\
\bar{w}_{yy}&=&\ez(1-\ez)t^{\ez-2}r^{2-2\sz}.
\end{eqnarray*}
We have
\begin{eqnarray*}
\mathrm{det}\;D^2_{r,y}\bar{w}&\ge&\ez^2r^{2-2\sz}t^{2\ez-2}[\sz(1-\ez)(3-\sz\ez)-(\sz\ez-2)^2]\\
&=&\ez^2r^{2-2\sz}t^{2\ez-2}[3\sz-4-\sz\ez(\sz-1)]
\end{eqnarray*}
and
$$\f{\bar{w}_r}{r}\ge\sz\ez t^{\ez}.$$
If we choose 
\begin{equation}\label{s2: *1}
\sz>\f{4}{3}
\end{equation}
and $\ez>0$ sufficiently small, then
\begin{eqnarray*}
\mathrm{det}\;D^2w_1(x',x_n)&\ge&c^{'n}C^{'2}c(n,\sz,\ez)r^{2-2\sz}t^{n\ez-2},\quad r:=|x'|,\;t:=C'x_n|x'|^{-\sz}.
\end{eqnarray*}
The right hand side of the last inequality is sufficiently large if we choose $C'>0$ sufficiently large.

Now we choose $c'$ small such that 
$$c'|x'|^{2}\le\mu|x'|^{2}\le u\quad\mathrm{on}\;\partial B^+_1\cap\{x_n\le\rz\}$$
and then $C'$ large such that
$$C'x_n|x'|^{-\sz}\ge 1\quad\mathrm{on}\;B^+_1\cap\{x_n=\rz\},$$
and
$$\mathrm{det}\;D^2w_1>\Lz\quad\mathrm{in}\;B^+_1\cap\{w_1>0\}.$$
Then we have $w_1\le u$ on $\partial(B^+_1\cap\{x_n\le\rz\})$ and $\mathrm{det}\;D^2w_1>\mathrm{det}\;D^2u$ on the set where $w_1>0$. Hence $u\ge w_1$ in $B^+_1\cap\{x_n\le\rz\}$.

It follows that
\begin{equation*}
S_h\cap\{x_n\le\rz\}\subset\l\{c'|x'|^{2}\le 2h\r\}\cup\l\{[1-(C'x_n|x'|^{-\sz})^{\ez}]\le\f{1}{2}\r\}
\end{equation*}
or
\begin{equation}\label{s2: *3}
S_h\cap\{x_n\le\rz\}\subset\{|x'|\le Ch^{\f{1}{2}}\}\cup\{x_n\ge c|x'|^{\sz}\}.
\end{equation}

\noindent$\mathbf{Step\;2:}$ Let $x^*_h$ be the center of mass of $S_h$ and denote $d_h=x^*_h\cdot e_n$, then we claim that for all small $h>0$ we have
\begin{equation}\label{s2: thm 13}
d_h\ge\tilde{c}h^{\f{\sz}{2}}
\end{equation}
for some small $\tilde{c}>0$.

Otherwise from \eqref{s2: *3} and Lemma \ref{Jo} we obtain
$$S_h\subset\{x_n\le C(n)\tilde{c}h^{\f{\sz}{2}}\le h^{\f{\sz}{2}}\}\cap\{|x'|\le C_1h^{\f{1}{2}}\}$$
for some large constant $C_1$. Then the function
$$w=\ez x_n+\l[\l|c_0\f{x'}{C_1h^{\f{1}{2}}}\r|^2+\l|C_0\f{x_n}{h^{\f{\sz}{2}}}\r|^2\r]\cdot h$$
is a lower barrier for $u$ in $S_h$ if $c_0$ is sufficiently small and $C_0$ is large.

Indeed, we have on $\partial S_h\cap\partial B^+_1\subset\{x_n=0\}$,
$$w=\l|c_0\f{x'}{C_1h^{\f{1}{2}}}\r|^{2}\cdot h\le\mu|x'|^{2}.$$
On $\partial S_h\cap B^+_1$,
$$w\le\ez+[c_0^{2}+C_0C(n)\tilde{c}]\cdot h\le h$$
if $c_0,\tilde{c}$ are sufficiently small.
Moreover, since $\sz>\f{4}{3}>1$ (see \eqref{s2: *1}), 
$$\mathrm{det}\;D^2w\ge 2^nh^n\l(\f{c_0}{C_1h^{\f{1}{2}}}\r)^{2(n-1)}\l(\f{C_0}{h^{\f{\sz}{2}}}\r)^2\ge\Lz$$
if $C_0$ is sufficiently large.

Hence, we obtain $w\le u$ in $S_h$. This is a contradiction with $\nabla u(0)=0$. Thus \eqref{s2: thm 13} is proved.

Since $S_h\cap\{x_n=0\}$ contains a ball in $\mathbb{R}^{n-1}$ of radius $(\mu h)^{\f{1}{1+\az}}$, we obtain from \eqref{s2: *3} and \eqref{s2: thm 13} that
$$\{te_n: 0\le t\le c_0h^{\f{\sz-1}{2}+\f{1}{1+\az}}\}\subset S_h$$
for some small constant $c_0>0$. Hence we only need to choose $\sz$ satisfying
\begin{equation}\label{s2: eq *4}
\f{\sz-1}{2}+\f{1}{1+\az}<1.
\end{equation}
By the assumption $\az>\f{1}{5}$, we can choose $\sz$ satisfying \eqref{s2: *1} and \eqref{s2: eq *4}. The conclusion of the theorem follows.

\section{Proof of Theorem \ref{thm 4}}\label{s6}

The proof of Theorem \ref{thm 4} is similar to that of Theorem \ref{thm 2}.\\

\noindent$\mathbf{Step\;1.}$ Let $y\in B^+_1$ and assume the maximal interior section $S_{\bar{h}}(y)$ is tangent to $\partial B^+_1$ at $0$. We prove that
\begin{equation}\label{thm 400}
S_{\bar{h}}(y)\subset B_{\bar{h}^{\ez_0}}
\end{equation}
for any $\bar{h}>0$ small, where $\ez_0$ is a small constant.

For any $h>0$ small, let $x^*_h$ be the center of mass of $S_h$ and $d_h:=x^*_h\cdot e_n$. Since $\overline{S_h}$ contains the convex set generated by $\{|x'|\le ch^{\f{1}{1+\az}},x_n=0\}$ and the point $x^*_h$, we can use the estimate of the upper bound of $|S_h|$ (see the proof of Theorem \ref{thm 2}) and obtain that
\begin{equation}\label{thm 40}
d_h\le C_0h^{\tau},\quad\quad\tau:=\f{n}{2}-\f{n-1}{1+\az}\in(0,\f{1}{2}]
\end{equation}
for some large constant $C_0>0$.

By Lemma \ref{Jo},
\begin{equation}\label{thm 41}
S_h\subset\{x_n\le C(n)C_0h^{\tau}\}\subset\{x_n\le C'h^{\tau}\}\quad\forall h>0.
\end{equation}
Choose $\sz=2>\f{4}{3}$ in \eqref{s2: *3} and we obtain that
\begin{equation}\label{thm 41'}
S_h\subset B_{Ch^{\f{\tau}{2}}}\quad\forall h>0.
\end{equation}

We have
$$S_{\bar{h}}(y)=\{x\in\Omega|\;v(x)<0\},$$
where 
$$v(x)=(u-l_0)(x)-[u_n(y)-u_n(0)]x_n.$$
Choose $h>0$ such that 
$$u_n(y)-u_n(0)=(2C')^{-1}h^{1-\tau},$$
where $C'$ is the constant in \eqref{thm 41}. Then we have by \eqref{thm 41} 
and \eqref{thm 41'} that
\begin{equation}\label{thm 41''}
S_{\bar{h}}(y)\subset S_h\subset B_{Ch^{\f{\tau}{2}}}.
\end{equation}
On the other hand, let $\tz>0$ be a small constant to be chosen below and denote 
$$\Lz_0:=\f{(1+\az')(1-\tau)}{\az'},$$
where $\az'\in(0,1)$ is the constant obtained by Theorem \ref{thm 3}. If $h\le 1$, by Theorem \ref{thm 3}, we can choose a point $z=te_n\in\partial S_{\tz h^{\Lz_0}}$ with $t\ge c(\tz h^{\Lz_0})^{\f{1}{1+\az'}}$. It follows that
\begin{eqnarray*}
v(z)&\le&\tz h^{\Lz_0}-(2C')^{-1}h^{1-\tau}c(\tz h^{\Lz_0})^{\f{1}{1+\az'}}\\
&=&\tz^{\f{1}{1+\az'}}h^{\Lz_0}[\tz^{\f{\az'}{1+\az'}}-(2C')^{-1}c]\\
&\le&-\tz^{\f{1}{1+\az'}}h^{\Lz_0}(4C')^{-1}c<0
\end{eqnarray*}
if $\tz>0$ is sufficiently small. On the other hand, if $h\ge 1$, then we can choose a point $z=te_n\in\partial S_{c'}$ with $t\ge cc'^{\f{1}{1+\az'}}$ and therefore
\begin{eqnarray*}
v(z)&\le&c'-(2C')^{-1}h^{1-\tau}cc'^{\f{1}{1+\az'}}\\
&=&c'^{\f{1}{1+\az'}}[c'^{\f{\az'}{1+\az'}}-(2C')^{-1}ch^{1-\tau}]\\
&\le&-c'^{\f{1}{1+\az'}}h^{1-\tau}(4C')^{-1}c<0
\end{eqnarray*}
if $c'>0$ is sufficiently small. Combining the last two estimates we obtain
\begin{equation}\label{thm 41'''}
\bar{h}\ge c|S_{\bar{h}}(y)|^{\f{2}{n}}\ge c|\min_{S_{\bar{h}}(y)}v|\ge c\min\{h^{\Lz_0},h^{1-\tau}\}.
\end{equation}
If $\bar{h}$ is sufficiently small, then \eqref{thm 41'''} together with \eqref{thm 41''} gives \eqref{thm 400} with $0<\ez_0<\f{\tau}{2\Lz_0}$.\\

\noindent$\mathbf{Step\;2.}$ Assume as in $\mathbf{Step\;1}$ that $S_{\bar{h}}(y)$ is a maximal interior section tangent to $\partial B^+_1$ at $0$. We prove that
\begin{equation}\label{thm 400**}
B_{\bar{h}^{1-\ez_1}}(y)\subset S_{\bar{h}}(y)
\end{equation}
for any $\bar{h}>0$ small, where $\ez_1>0$ is a small constant.

By Lemma \ref{Jo}, $S_{\bar{h}}(y)$ is equivalent to an ellipsoid $E$ centered at $y$ of axes $0<\mu_1\le\mu_2\le\dots\le\mu_n$.  We only need to prove that
\begin{equation}\label{thm 42}
\mu_1\ge\bar{h}^{1-\ez_1}.
\end{equation}
Assume by contradiction that 
$$\mu_1<\bar{h}^{1-\ez_1}.$$
Then as in the proof of Theorem \ref{thm 2} we have
\begin{equation}\label{thm 43}
v\ge c\bar{h}^{\ez_1}|(x-y)\cdot\nu|-2\bar{h}\quad\mathrm{in}\;B^+_1,
\end{equation}
where $\nu$ is a unit vector.

We have 
\begin{equation}\label{thm 44'}
v|_{\partial B^+_1\cap\{x_n=0\}}=\varphi-\varphi(0)-\nabla_{x'}\varphi(0)\cdot x'\le C|x'|^{1+\az}.
\end{equation}

We only need to consider the cases: $|\nu'|\ge\bar{h}^{\f{\ez_0\az}{4}}$ and $|\nu'|\le\bar{h}^{\f{\ez_0\az}{4}}$, where $\ez_0$ is the constant in  \eqref{thm 400}. 

\noindent$\mathbf{Case\;1.}$ $|\nu'|\ge\bar{h}^{\f{\ez_0\az}{4}}$. We choose $x=(x',0)$ with $x'=\bar{h}^{\f{\ez_0}{2}}\f{\nu'}{|\nu'|}$. Then we have
$$|y|\le\bar{h}^{\ez_0}\le\bar{h}^{\f{\ez_0}{2}}|x'|,$$
which implies 
\begin{eqnarray*}
\bar{h}^{\ez_1}(x-y)\cdot\nu&\ge&\bar{h}^{\ez_1}\l(x'\cdot\nu'-|y|\r)\\
&\ge&\bar{h}^{\ez_1+\f{\ez_0\az}{4}}(1-C\bar{h}^{\f{\ez_0}{2}-\f{\ez_0\az}{4}})|x'|\\
&\ge&\f{\bar{h}^{\ez_1+\f{\ez_0\az}{4}}}{2}|x'|
\end{eqnarray*}
if $\ez_1,\ez_0>0$ are sufficiently small and $\bar{h}$ is small. Choose $0<\ez_1<\f{\ez_0\az}{4}$, then by \eqref{thm 43} and \eqref{thm 44'}, we reach a contradiction. Thus \eqref{thm 42} holds.

\noindent$\mathbf{Case\;2.}$ $|\nu'|\le\bar{h}^{\f{\ez_0\az}{4}}$. 

By \eqref{thm 41'''} we find that $h$ is bounded above, hence we have $v+\bar{h}\ge 0$ in $B^+_1$ and
$$v\ge c|x'|^2\quad\mathrm{on}\;\partial B^+_1\cap\{x_n\le c_*\}$$
if $c_*>0$ is sufficiently small, then similar to the proof of Theorem \ref{thm 3}, the function
$$\tilde{w}_1:=\tilde{c}'|x'|^2[1-(\tilde{C}'x_n|x'|^{-2})^{\ez}]^+$$
is a lower barrier for both $u-l_0$ and $v+\bar{h}$ in $B^+_1\cap\{x_n\le c_*\}$, where $\ez,\tilde{c}'$ are some small constants and $\tilde{C}'$ is a large constant. 

On $\{x_n=\f{1}{2\tilde{C}'}|x'|^2\}\cap\{|x'|\le\f{1}{2}\}$ we have
$$u-l_0\ge\tilde{w}_1\ge(1-2^{-\ez})\tilde{c}'|x'|^2,$$
it follows that
\begin{eqnarray}\label{thm 45}
v&\ge&(1-2^{\ez})\tilde{c}'|x'|^2-(2C')^{-1}h^{1-\tau}\f{1}{2\tilde{C}'}|x'|^2\nonumber\\
&\ge&\f{1}{2}(1-2^{\ez})\tilde{c}'|x'|^2
\end{eqnarray}
if $\tilde{C}'$ is sufficiently large.

On $\{\bar{h}^{\f{\ez_0\az}{4}}=x_n>\f{1}{2\tilde{C}'}|x'|^2\}$, we have by \eqref{thm 43}
\begin{eqnarray}\label{thm 46}
v&\ge&c\bar{h}^{\ez_1}\l(x_n|\nu\cdot e_n|-\bar{h}^{\f{\ez_0\az}{4}}|x'|-|y|\r)-2\bar{h}\nonumber\\
&\ge&c\bar{h}^{\ez_1}\l(\f{1}{2}x_n-\bar{h}^{\f{\ez_0\az}{4}}(2\tilde{C}'x_n)^{\f{1}{2}}-\bar{h}^{\ez_0}\r)-2\bar{h}\nonumber\\
&\ge&c\bar{h}^{\ez_1}x_n.
\end{eqnarray}

Define 
$$w=\dz x_n+c_2\bar{h}^{\ez_1}|x'|^2+\f{C_2}{\bar{h}^{\ez_1(n-1)}}x_n^2,$$
then $w$ is a lower barrier for $v$ in $\{\f{1}{2\tilde{C}'}|x'|^2\le x_n\le\bar{h}^{\f{\ez_0\az}{4}}\}$ if $\dz,c_2>0$ are sufficiently small and $C_2$ is large. Indeed, on $\{\bar{h}^{\f{\ez_0\az}{4}}\ge x_n=\f{1}{2\tilde{C}'}|x'|^2\}$, we have by \eqref{thm 45}
$$w\le 2\dz\tilde{C}'|x'|^2+c_2\bar{h}^{\ez_1}|x'|^2+C\bar{h}^{\f{\ez_0\az}{4}-\ez_1(n-1)}|x'|^2\le v$$
if $\dz,c_2$ are small and $\ez_1(n-1)<\f{\ez_0\az}{4}$.

On $\{\bar{h}^{\f{\ez_0\az}{4}}=x_n>\f{1}{2\tilde{C}'}|x'|^2\}$, we have by \eqref{thm 46}
$$w\le\dz x_n+c_2\bar{h}^{\ez_1}2\tilde{C}'x_n+C_2\bar{h}^{\f{\ez_0\az}{4}-\ez_1(n-1)}x_n\le v$$
if $\dz,c_2$ are small and $\ez_1n<\f{\ez_0\az}{4}$.

Therefore, $v\ge w\ge\dz x_n$ in $\{\f{1}{2\tilde{C}'}|x'|^2\le x_n\le\bar{h}^{\f{\ez_0\az}{4}}\}$. Since 
$$0=v(y)+\bar{h}\ge\tilde{w}_1(y),$$
we have
$$\bar{h}^{\ez_0}\ge y_n\ge\f{1}{\tilde{C}'}|y'|^{2}.$$
Hence we reach a contradiction since $v(y)<0$. Thus \eqref{thm 42} holds. Combining these two cases, we obtain \eqref{thm 400**}.\\

\noindent$\mathbf{Step\;3.}$ Let $c_0>0$ be a small constant to be chosen below and $y\in B^+_{c_0}$. Assume the maximal interior section $S_{\bar{h}}(y)$ is tangent to $\partial B^+_1$ at $y_0\in\partial B^+_1$. We prove that
\begin{equation}\label{thm 480}
B_{\bar{h}^{1-\ez_1}}(y)\subset S_{\bar{h}}(y)\subset B_{\bar{h}^{\ez_0}}(y_0),
\end{equation}
where $\ez_1,\ez_0$ are the constants in \eqref{thm 400} and \eqref{thm 400**}.

We only need to prove that
\begin{equation}\label{thm 48}
y_0\in\{|x'|\le\f{1}{2},x_n=0\}.
\end{equation}
Then it follows from $\mathbf{Step\;1}$-$\mathbf{2}$ (with $0$ replaced by $y_0$) that \eqref{thm 480} holds.

Now we prove \eqref{thm 48}.

Since $S_{\bar{h}}(y)$ is equivalent to an ellipsoid centered at $y$, we obtain
\begin{equation*}
y_0\cdot e_n\le Cy_n\le Cc_0.
\end{equation*}

Assume in contradiction that \eqref{thm 48} does not hold, then we have 
$$y_0\in\{\f{1}{2}\le|x'|\le 1,x_n=0\}\cup\{|x|=1,0<x_n\le Cc_0\}.$$
Recall that $u$ separates quadratically from $l_0$ on $\partial B^+_1$ in a neighborhood of $0$. Hence if $c_0$ is sufficiently small, then \begin{equation}\label{thm 47'}
(u-l_0)(y_0)\ge c|y'_0|^2\ge c_1
\end{equation}
for some $c_1$ small. Since
\begin{eqnarray*} 
0\le (u-l_y)(x)=(u-l_0)(x)+(u-l_y)(0)+(\nabla u(0)-\nabla u(y))\cdot x
\end{eqnarray*}
and $u$ is Lipschitz continuous, we have on $\partial B^+_1\cap\{x_n=0\}$,
$$[\nabla_{x'}u(y)-\nabla_{x'}\varphi(0)]\cdot x'\le C[|x'|^{1+\az}+|y|].$$
Choose $x'=\f{\nabla_{x'}u(y)-\nabla_{x'}\varphi(0)}{|\nabla_{x'}u(y)-\nabla_{x'}\varphi(0)|}|y|^{\f{1}{2}}$, we obtain that
\begin{equation}\label{thm 47''}
|\nabla_{x'}u(y)-\nabla_{x'}\varphi(0)|\le C|y|^{\f{\az}{2}}.
\end{equation}
The convexity of $u$ implies that $u_n(0)$ is bounded above. Since $u$ is Lipschitz continuous, $|\nabla u(y)|$ is bounded and moreover, for any $x\in B^+_1$,
$$u(x)=u(x)-\varphi(x')+\varphi(x')\ge-Cx_n+\varphi(0)+\nabla_{x'}\varphi(0)\cdot x',$$
which implies that $u_n(0)$ is bounded below.

We have
\begin{eqnarray*} 
v&=&u-l_y-(u-l_y)(y_0)\\
&=&u-l_0-(u-l_0)(y_0)-(\nabla u(y)-\nabla u(0))\cdot (x-y_0).
\end{eqnarray*}
On $\partial B^+_1\cap\{x_n=0\}$, since $v\ge 0$, we have by \eqref{thm 47'} and \eqref{thm 47''}
\begin{eqnarray}\label{*}
C|x'|^{1+\az}&\ge&\varphi-\varphi(0)-\nabla_{x'}\varphi(0)\cdot x'\nonumber\\
&\ge&(u-l_0)(y_0)+(\nabla_{x'}u(y)-\nabla_{x'}\varphi(0))\cdot (x'-y'_0)\nonumber\\
&+&(u_n(y)-u_n(0))(-y_0\cdot e_n)\nonumber\\
&\ge&c_1-Cc_0^{\f{\az}{2}}-Cc_0\ge\f{c_1}{2}
\end{eqnarray}
if $c_0$ is sufficiently small. We reach a contradiction. Hence \eqref{thm 48} holds. 

Using \eqref{thm 480} and similar arguments as in $\mathbf{Step\;3}$ in the proof of Theorem \ref{thm 2} we can prove
$$[\nabla u]_{C^{\bz}(\overline{B^+_{c_0}})}\le C$$
for some $\bz\in(0,1)$ depending only on $n,\lz,\Lz$ and $\az$.

\end{document}